\documentclass[a4paper,11pt,reqno]{amsart}

\usepackage{amssymb}

\def\Z{{\mathbb Z}}

\def\subo{{_{\bar 0}}}
\def\subuno{{_{\bar 1}}}

\DeclareMathOperator{\rank}{rank}

\newtheorem{defi}{\bf Definition}[section]
\newtheorem{theo}[defi]{\bf Theorem}

\newtheorem{coro}[defi]{\bf Corollary}

\newtheorem{lem}[defi]{\bf Lemma}

\def\hregla{\hrule height.1pt}

\def\hreglon{\hrule height1pt}
\def\vreglon{\vrule height 12pt width1pt depth 4pt}

\def\vregla{\vrule width.1pt}

\begin{document}

\title[The Kac Jordan superalgebra]{The Kac Jordan superalgebra: automorphisms and maximal subalgebras}

\author{Alberto Elduque}
 \address{Departamento de Matem\'aticas, Universidad de
Zaragoza, 50009 Zaragoza, Spain}
 \email{elduque@unizar.es}

\author{Jes\'us Laliena}
 \address{Departamento de Matem\'aticas y Computaci\'on,
 Universidad de La Rioja,
 26004, Logro\~no, Spain}
 \email{jesus.laliena@dmc.unirioja.es}

\author{Sara Sacrist\'an}
 \address{Departamento de Matem\'aticas y Computaci\'on,
 Universidad de La Rioja,
 26004, Logro\~no, Spain}

\thanks{The first and second authors
have been supported by the Spanish Ministerio de
 Educaci\'{o}n y Ciencia
 and FEDER (MTM 2004-081159-C04-02). The first author also acknowledges
 support of the Diputaci\'on General de Arag\'on (Grupo de
 investigaci\'on de \'Algebra)}




\begin{abstract}

In this note the group of automorphisms of the Kac Jordan
superalgebra is described, and used to classify the maximal
subalgebras.

\end{abstract}

\maketitle

\section{Introduction.}


Finite dimensional simple Jordan superalgebras over an algebraically
closed field of characteristic zero were classified by V.~Kac in
1977 \cite {Kac} (see also  Kantor  \cite {Kan}, where a missing
case is added). Among these superalgebras we find the ten
dimensional Kac Jordan superalgebra, $K_{10}$, which is exceptional
 (see
\cite{McC2} and \cite{Sht}) and plays a significant role (see
\cite{MedZel}).

Here we are interested in describing the group of automorphisms of $K_{10}$ and, as a consequence,
in classifying the maximal subalgebras of
the Kac superalgebra over an algebraically closed field of characteristic zero.

In two previous works (see \cite{mio}, \cite{otro}), the authors
have given a description of maximal subalgebras of finite
dimensional central simple superalgebras which are either
associative or associative with superinvolution and also a
description of maximal subalgebras of finite dimensional special
simple Jordan superalgebras  over an algebraically closed field of
characteristic zero. These papers are thus in the spirit of previous
work by  E.~Dynkin (\cite{Dy1}, \cite{Dy2}), M.~Racine (\cite{Ra1},
\cite{Ra2}) or A.~Elduque (\cite{E}) about maximal subalgebras of
different classes of algebras.

Let $F$ be a field of characteristic not two. This assumption will
be kept throughout the paper. Recall that a superalgebra is a
$\Z_2$-graded algebra $A= A\subo \oplus A\subuno$ ($A_\alpha A_\beta
\subseteq A_{\alpha + \beta}\ \forall\alpha, \beta \in \Z_2$). If
$a\in A_\alpha$ we say that $a$ is a \emph{homogeneous element} and
we use the notation $\bar a = \alpha$. Elements belonging to
$A_{\bar0}$ are called \emph{even} elements and the ones in
$A_{\bar1}$ \emph{odd} elements. A superalgebra $A$ is said to be
non trivial if $A\subuno \not = 0$, and it is called \emph{simple}
if it contains no proper nontrivial graded ideals  and $A^2\neq 0.$

\bigskip

A Jordan algebra is an algebra satisfying the following identities:
\begin{gather}
 xy=yx, \\
(x^2y)x=x^2(yx).
\end{gather}
The last identity can be written as $(x^2,y,x)=0$, where $(x,y,z)$
denotes the associator $(xy)z-x(yz)$. Following the standard
procedure, a superalgebra is a \emph{Jordan superalgebra} if its
Grassmann envelope is a Jordan algebra. In particular, a Jordan
superalgebra satisfies the identities
\begin{gather}
xy=(-1)^{\bar x\bar y}yx\\
\sum_{\substack{\textrm{cyclic}\\ x,y,t}}(-1)^{(\bar x+\bar z)\bar
t} (xy,z,t)=0
\end{gather}
for homogeneous elements $x,y,z,t$. (If the characteristic is $\ne
2,3$, these two identities characterize Jordan superalgebras.)

%
%

\medskip

The even part of a Jordan superalgebra is a Jordan algebra, while
the odd part is a Jordan bimodule for the even part.


Associative superalgebras are just $\Z_2$-graded associative
algebras, but note that Jordan superalgebras are not $\Z_2$-graded
Jordan algebras.

The following examples of Jordan algebras will be useful in the
sequel:

\begin{enumerate}
\item[{\rm i)}] Let $A$ be an associative algebra over a ground field of characteristic $\not= 2$.
The new operation $a \cdot b = \frac{1}{2} (ab+ba)$ defines a
structure of a Jordan algebra on $A$. It is denoted by $A^+$.

\smallskip

\item[{\rm ii)}] Let $V$ be a vector space over $F$ with a nondegenerate symmetric bilinear form $(
\ , \ ) \colon V\times V \to F$. The direct sum $F1+V$ with the
product  $(\lambda 1 + v) (\mu 1+ w)= (\lambda \mu + (v,w))1+
(\lambda w + \mu v)$ is a Jordan algebra, and it is called the
Jordan algebra of a nondegenerate bilinear form.

\end{enumerate}

\smallskip

 Simple finite dimensional Jordan algebras over an algebraically closed field were classified by
 P.~Jordan, J.~von Neumann, E.~Wigner (\cite {JNW}) and by A.~Albert
(\cite {A}). The examples given above are two of the four types of
algebras given in that classification.

\bigskip

Also the following examples of Jordan superalgebras will be needed
later on:

\begin{enumerate}
\item[{\rm 1)}] $J=K_3,$ the \emph{Kaplansky superalgebra}:

 $J_{\bar{0}}=Fe,\quad J_{\bar{1}}=Fx+Fy,$

  $e^2=e, \quad e \cdot x=\frac{1}{2} x, \quad  e \cdot y = \frac{1}{2}y, \quad x \cdot y=e.$

\bigskip

\item[{\rm 2)}]  The family of superalgebras $J=D_t,$ with $t\in F
\setminus\{0\}$:

 $J_{\bar{0}}=Fe +Ff,$\quad  $J_{\bar{1}} = Fu+Fv$

 $e^2 =e, \quad f^2=f, \quad  e \cdot f=0,\quad  e \cdot u=\frac{1}{2}u, \quad e \cdot v=\frac{1}{2}v, \quad
f \cdot u=\frac{1}{2}u,$

$ f \cdot v=\frac{1}{2}v, \quad u \cdot v=e+tf.$

\bigskip

\item[{\rm 3)}] $J=K_{10},$ the \emph{Kac superalgebra}, which  will be
described in more detail in the following section.

\bigskip

\item [{\rm 4)}] The {\it superalgebra of a superform}: Let $V= V_{\bar{0}}\oplus V_{\bar{1}}$ be a
graded vector space over $F$, and let $( \ ,\ )$ be a nondegenerate
supersymmetric bilinear superform on $V$, that is, a nondegenerate
bilinear map which is symmetric on $V_{\bar 0}$, skewsymmetric on
$V_{\bar 1}$, and $V_{\bar 0}, V_{\bar 1}$ are orthogonal relative
to $(\ , \ )$. Now consider $J_{\bar 0} =Fe+V_{\bar 0}$, $J_{\bar
1}=V_{\bar 1}$ with $e\cdot x=x$, $v\cdot w = (v,w) e$ for every
$x\in J$, $v,w \in V$. This superalgebra, $J$, is called  the
superalgebra of a superform.

\bigskip

\item [{\rm 5)}] $A^{+},$ with $A$ an
associative superalgebra over a field of characteristic not 2, where
the product operation in $A$ has been changed to $a_i\cdot b_j=
\frac{1}{2}(a_i b_j + (-1)^{ij}b_j a_i).$

\end{enumerate}

\bigskip

As we have mentioned at the beginning of this Introduction, V. Kac
and I. Kantor gave the classification of nontrivial simple Jordan
superalgebras. The four examples given above are also examples of
simple Jordan superalgebras, assuming in 3) that the characteristic
is not 3, and in 4) that $A$ is a simple associative superalgebra
(see \cite{Wall}).

\bigskip

In this paper, our purpose is to describe the group of automorphisms
of $K_{10}$ and also to classify its maximal subalgebras.

\bigskip

\section{Kac superalgebra}
The Kac superalgebra, $K_{10}= J_{\bar{0}} \oplus J_{\bar{1}},$ is a
Jordan superalgebra with 6-dimensional even part  and 4-dimensional
odd part:
$$
J_{\bar{0}}= (Fe + Fa+ Fb + \sum_{i=1}^2 Fc_i) \oplus Ff,
\hspace{1cm} J_{\bar{1}}=  \sum_{i=1}^2 (Fp_i+ Fq_i)
$$
\noindent and product given in Table \ref{ta:K10}.

\begin{table}[h!]
{\footnotesize
$$ \vbox{\offinterlineskip
\halign{\hfil$#$\enspace\hfil\vreglon
 &\hfil\enspace$#$\enspace\hfil
 &\hfil\enspace$#$\enspace\hfil
 &\hfil\enspace$#$\enspace\hfil
 &\hfil\enspace$#$\enspace\hfil
 &\hfil\enspace$#$\enspace\hfil
 &\hfil\enspace$#$\quad\hfil\vregla
 &\hfil\qquad\enspace$#$\quad\enspace\hfil
 &\hfil\quad\enspace$#$\quad\enspace\hfil
 &\hfil\quad\enspace$#$\quad\enspace\hfil
 &\hfil\quad\enspace$#$\quad\enspace\hfil\cr
 \omit\hfil\vrule width1pt depth 4pt
   &e&a&b&c_1&c_2&f&p_1&p_2&q_1&q_2\cr
 \noalign{\hreglon}
 e&e&a&b&c_1&c_2&0&
   \tfrac{1}{2}p_1&\tfrac{1}{2}p_2&\tfrac{1}{2}q_1&\tfrac{1}{2}q_2\cr
 a&a&4e&0&0&0&0&p_1&p_2&-q_1&-q_2\cr
 b&b&0&-4e&0&0&0&q_1&q_2&-p_1&-p_2\cr
 c_1&c_1&0&0&0&2e&0&0&q_1&0&p_1\cr
 c_2&c_2&0&0&2e&0&0&q_2&0&p_2&0\cr
 \omit $f$\enspace\hfil\vrule height 12pt width1pt depth 6pt& 0&0&0&0&0&f&\tfrac{1}{2}p_1&
   \tfrac{1}{2}p_2&\tfrac{1}{2}q_1&\tfrac{1}{2}q_2\cr
 \noalign{\hregla}
 p_1&\tfrac{1}{2}p_1&p_1&q_1&0&q_2&\tfrac{1}{2}p_1&
   0&\hidewidth a+2(e-3f)\hidewidth&2c_1&b\cr
 p_2&\tfrac{1}{2}p_2&p_2&q_2&q_1&0&\tfrac{1}{2}p_2&
   \hidewidth -a-2(e-3f)\hidewidth&0&-b&-2c_2\cr
 q_1&\tfrac{1}{2}q_1&-q_1&-p_1&0&p_2&\tfrac{1}{2}q_1&
   -2c_1&b&0&\hidewidth a-2(e-3f)\hidewidth\cr
 q_2&\tfrac{1}{2}q_2&-q_2&-p_2&p_1&0&\tfrac{1}{2}q_2&
   -b&2c_2&\hidewidth -a+2(e-3f)\hidewidth &0\cr
}}
$$}
\smallskip
\caption{$K_{10}$}\label{ta:K10}
\end{table}

\bigskip

This basis of $K_{10}$ is obtained by using the one in
\cite[p.~378]{RaZ} with $a=v_1+v_2$, $b=v_1-v_2$, $c_1=v_3$,
$c_2=v_4$, $p_1=x_1-y_2$, $p_2=x_2+y_1$, $q_1=x_1+y_2$,
$q_2=x_2-y_1$.

\medskip

 From \cite {Ki} we know that $J_{\bar{0}}= J(V,Q)\oplus Ff,$ where
$J(V,Q)$ is the Jordan algebra of a nondegenerate bilinear form $Q,$
$e$ is the identity element in $J(V,Q),$ and $V$ is the vector space
with basis $\{ a, b, c_1, c_2\}.$ We will use this presentation in
Section 4.

\vspace{1cm}

In \cite {N-B-E} G.~Benkart and A.~Elduque gave a realization of
$K_{10}$ which  enables to check directly that it is indeed a Jordan
superalgebra. We introduce this construction that will be very
useful for us in Section 3.

First we pick up the Kaplansky superalgebra, $K_3,$ over the field
$F$ and we define the following supersymmetric bilinear form (that
is, the even and odd part are orthogonal and the form is symmetric
in the even part and alternating in the odd part):
$$ (e|e)= \frac{1}{2}, \quad (x|y)=1, \quad ((K_3)_{\bar{0}}|(K_3)_{\bar{1}})=0.$$

\noindent Consider now the vector space over $F$: $F\cdot 1 \oplus
(K_3 \otimes _F K_3)$ and  define on it the product:
$$ (a\otimes b)(c\otimes d)= (-1)^{\bar{b}\bar{c}} (ac\otimes bd - \frac{3}{4} (a|c)(b|d)1)$$

\noindent with $a,b,c,d\in K_3$ homogeneous elements, and where $1$
is a formal identity element. Then,

\bigskip

\begin{theo}\label{th:BeEl} (Benkart, Elduque)
$K_{10}$ is isomorphic to $F \cdot 1 \oplus (K_3\otimes K_3)$ by
means of the linear map given by:
\begin{align*}
 e & \mapsto  \frac{3}{2} \cdot 1 - 2 e\otimes e,&
 f & \mapsto  -\frac{1}{2} \cdot 1 + 2 e\otimes e,\\
a  & \mapsto  -4x \otimes x-y\otimes y,&
 b &\mapsto  -4x\otimes x + y\otimes y,\\
 c_1 &\mapsto  2x\otimes y,&
c_2 & \mapsto  -2y\otimes x,\\
p_1 & \mapsto  4x\otimes e - 2e\otimes y,&
p_2 &\mapsto -4e\otimes x -2y\otimes e,\\
q_1 &\mapsto 4x\otimes e + 2e\otimes y,&
 q_2 &\mapsto -4e\otimes x+ 2y\otimes e.
\end{align*}

\end{theo}

\bigskip

If we denote by $W$ the vector space $(K_3)_{\bar 1}$ we notice that
$V= W \otimes W$ is endowed with the bilinear nondegenerate
symmetric form given by $b(s\otimes t, u\otimes v)= (s|u) (t|v)$. We
note also that the corresponding decomposition of $(K_{10})_{\bar
0}$ into the direct sum of two simple ideals according to the one
given above is $(K_{10})_{\bar 0}=I\oplus J$, where $I=F\cdot(-1/2+2
e\otimes e)$ and $J=F\cdot (3/2-2 e\otimes e) \oplus V$.

\bigskip

\section {Automorphisms of $K_{10}$}

Here by automorphism we mean automorphism of graded algebras.

It is easy to compute the group of automorphisms of the examples 1),
2) and 4) of Jordan superalgebras given in Section 1. So, with the
notation of Section 2, $Aut(K_3)=Sp(W)$, the symplectic group (or
also $Aut(K_3)=SL(W)$, the special linear group). And
$Aut(D_t)=Sp(U)$ where $U$ is the vector space generated by
$\{u,v\}$. If $J$ is the superalgebra of a superform then $Aut(J) =
O(n,F)\oplus Sp(m,F)$ where $O(n,F)$ denotes the orthogonal group of
$V_{\bar 0}$, with $\dim V_{\bar 0}=n$, related to $(\ , \
)|_{V_{\bar 0}\times V_{\bar 0}}$, and $Sp(m,F)$ denotes the
symplectic group of $V_{\bar 1}$, with $\dim V_{\bar 1}=m$, related
to $(\ , \ )|_{V_{\bar 1}\times V_{\bar 1}}$.

\bigskip

In what follows in this section it will assumed that $F$ is a field
of characteristic $\not= 2,3$  such that $F^2=F$.

\bigskip

Let $C_2 = \{1, \epsilon\}$ be the cyclic group of order $2$
($\epsilon^2 =1$) and consider the wreath product $G = Sp(W)\wr
 C_2$ (that is, $G$ is the
semidirect product $(Sp(W) \times Sp(W)) \rtimes C_2$ and $(f,g)
\epsilon = \epsilon (g,f)$ for any $f,g \in Sp(W)$).

The following maps are clearly group homomorphisms:

\begin{eqnarray*}
\Psi:& Aut(K_{10}) &\longrightarrow O(V,b)\\
& \varphi &\longmapsto \varphi|V \\
\\
\tilde{\Psi}:& G &\longrightarrow O(V,b)\\
& (f, g) &\longmapsto f \otimes g: W \otimes W \rightarrow \hspace{0.25cm} W \otimes W \\
&        &    \hspace{2.25cm} s \otimes t \hspace{0.3cm}\mapsto
f(s) \otimes g(t)\\
& \epsilon      & \longmapsto \hspace{0.5cm}\hat{\epsilon}
\hspace{0.25cm}: W \otimes W \rightarrow
\hspace{0.25cm} W \otimes W \\
&  &   \hspace{2.25cm} s \otimes t \hspace{0.3cm}\mapsto
\hspace{0.25cm}-t \otimes x\\
\\
\Phi:& G &\longrightarrow Aut(K_{10})\\
& (f, g) &\longmapsto \hspace{0.25cm} \Phi_{(f, g)}:
\hspace{0.25cm}\left\{\begin{array}{rl} 1 &
\rightarrow 1 \\
a \otimes b & \mapsto \tilde{f}(a) \otimes \tilde{g}(b)
\end{array}\right.
\\
& \epsilon &\longmapsto \hspace{0.5cm} \delta :
\hspace{0.75cm}\left\{\begin{array}{rl} 1 &
\rightarrow 1 \\
a \otimes b & \mapsto (-1)^{\bar{a}\bar{b}} b \otimes a
\end{array}\right.
\\
\end{eqnarray*}

\noindent where for any $f \in Sp(W)$, $\tilde{f}$ denotes the
automorphism of $K_3$ such that $\tilde{f}(e)= e,$ $\tilde{f}(w)=
f(w)$ $\forall w \in W.$

Note that $\Psi \circ \Phi = \tilde{\Psi}.$

\medskip

\begin{lem}
$\tilde{\Psi}$ is onto with $Ker \tilde{\Psi} = \{ \pm (id,
id)\}.$
\end{lem}
\begin{proof}[Proof:]
Let $\{u,v\}$ be a symplectic basis of $W$ (so $(u \thickspace \vert
\thickspace v)= 1),$ then if $\alpha = u \otimes s + v \otimes t$ is
an isotropic vector in $V$ (that is, $b(\alpha, \alpha)=0)$, then $0
= b(\alpha,\alpha) = (s \thickspace \vert \thickspace t)$, so $s, t$
are linearly dependent. This shows that
$$\{ \hbox{isotropic vectors of }  V\} = \{ s \otimes t \thickspace \vert \thickspace s, t \in W\}.$$

Now let us show that $\tilde{\Psi}$ is onto. For any $\varphi \in
O(V, b),$ $\varphi(u \otimes u)$ is isotropic (because so is $u
\otimes u$), so there are $s, t \in W$ with $\varphi(u \otimes u)= s
\otimes t \neq 0.$

Take $f, g \in Sp(W)$ with $f(s) = u,$ $g(t) = u.$ Hence $\varphi
\in im \tilde{\Psi}$ if and only if
$$
(f \otimes g) \circ \varphi = \tilde \Psi(f, g) \circ \varphi \in im \tilde \Psi
$$
and, therefore, it can be assumed that $\varphi(u \otimes u)= u
\otimes u.$

\bigskip

Now $\varphi (v \otimes v) = s \otimes t$ for some $s, t \in W$
with
$$1 = b (u \otimes u, v \otimes v)= b (\varphi(u \otimes u), \varphi(v \otimes v))=
b (u \otimes u, s \otimes t)= (u \thickspace \vert \thickspace s)(u
\thickspace \vert \thickspace t)$$ \noindent and there exists $0
\neq \alpha \in F$ such that $(u \thickspace \vert \thickspace s)=
\alpha,$  $(u \thickspace \vert \thickspace t)= \alpha^{-1}$ or $(u
\thickspace \vert \thickspace \alpha^{-1}s)=1,$ $(u \thickspace
\vert \thickspace \alpha t)= 1.$ Thus there are $f, g \in Sp(W)$
such that $f(u) = u,$ $f(v) = \alpha^{-1}s, $ $g(u) = u,$ $g(v) =
\alpha t $, and changing $\varphi$ to $(f^{-1} \otimes g^{-1}) \circ
\varphi$ we may assume that $\varphi(u \otimes u) = u \otimes u,$
$\varphi(v \otimes v) = v \otimes v.$

Then $\varphi(u \otimes v)$ is isotropic and orthogonal to both $u
\otimes u$ and $v \otimes v,$ so it is a scalar multiple of either
$u \otimes v$ or $v \otimes u.$  By using $\tilde \Psi (\epsilon)
= \hat \epsilon$ we may assume that $\varphi (u \otimes v)= \gamma
u \otimes v$ for some $0 \neq \gamma \in F$ and then, necessarily
$\varphi (v \otimes u)= \gamma^{-1} v \otimes u.$ Let $\mu \in F =
F^2$ with $\mu^2 = \gamma,$ then $\varphi = f \otimes g$ where
$f(u) = \mu u,$ $f(v) = \mu^{-1} v,$ $g(u) = \mu^{-1} u,$ $g(v) =
\mu v,$ so that $\varphi \in im \tilde{\Psi}.$

The assertion about the kernel is clear.
\end{proof}

\bigskip

\begin{lem}
$\Psi$ is onto with $Ker \Psi = \{id, \tau\},$ where $\tau$ is the
grading automorphism $(\tau (z) = (-1)^{\bar z} z)$ for any
homogeneous $z \in K_{10})$.
\end{lem}
\begin{proof}[Proof:] Since $\Psi \circ \Phi = \tilde{\Psi}$ and
$\tilde{\Psi}$ is onto, so is $\Psi.$ Finally, if $\varphi \in Aut
K_{10}$ belongs to the kernel of $\Psi,$ then $\varphi |
(K_{10})_{\bar{0}} = Id,$ so $\varphi (z_{\bar{0}} \cdot
z_{\bar{1}}) = z_{\bar{0}} \cdot \varphi(z_{\bar{1}})$ for any
$z_{\bar{0}} \in (K_{10})_{\bar{0}}$ and $z_{\bar{1}} \in
(K_{10})_{\bar{1}}.$

For any $s \in W,$ $s \otimes e$ is, up to scalars, the unique odd
element annihilated by $s \otimes W$ ($\subseteq
(K_{10})_{\bar{0}}$), hence $\varphi (s \otimes e) = \beta s
\otimes e$ for some $0 \neq \beta \in F.$

Take $t \in W$ with $(s \thickspace \vert \thickspace t) = 1,$
then $\forall z \in W$ $(t \otimes z) \cdot (s \otimes e)=
\frac{1}{2} e \otimes z$ so $\varphi (e \otimes z)= 2 (t \otimes
z) \cdot \varphi(s \otimes e) = \beta (e \otimes z)$ for any $z
\in W.$ It follows that $\varphi \vert (K_{10})_{\bar{1}} =
\varphi \vert (e \otimes W \oplus W \otimes e) = \beta Id$ and,
since $\varphi$ is an automorphism, $\beta^2 = 1$ so $\varphi$ is
either the identity or the grading automorphism $\tau.$
\end{proof}

\bigskip

\begin{theo}
$\Phi$ is an isomorphism.
\end{theo}
\begin{proof}[Proof:]
Since $\Phi (- id, -id)= \tau ,$  $\Phi(Ker \tilde \Psi) = Ker \Psi.$
Now the assertion follows from the fact that $\Psi \circ \Phi =
\tilde \Psi,$ together with Lemmata 1 and 2.
\end{proof}

\bigskip

\section{Maximal subalgebras of $K_{10}$}


In what follows the word subalgebra will be used in the graded
sense, so any sub\-al\-ge\-bra is graded. Also in this paragraph we
consider that $K_{10}$ is a superalgebra over $F$, an algebraically
closed field of characteristic $\ne 2,3$.

\bigskip

\begin{theo}

Any maximal subalgebra of $K_{10}$ is, up to an automorphism of
$K_{10},$ one of the following:

\begin{enumerate}

\item[{\rm (i)}] $(K_{10})_{\bar{0}}.$

\item[{\rm (ii)}] The subalgebra with basis $\{ e,f,a,p_1,p_2\}.$

\item[{\rm (iii)}] The subalgebra with basis
$\{e,f,a+b,c_1,p_1,q_1, p_2+q_2\}.$

\item[{\rm (iv)}] The subalgebra with basis
$\{e,f,a,b,c_1,p_1,q_1\}.$
\end{enumerate}

\end{theo}

\begin{proof}[Proof:]  Note that if $B$ is a maximal subalgebra of
$K_{10}$, $1=e+f\in B$, as $F1+B$ is a subalgebra of $K_{10}$ and
$B$ is an ideal in $F1+B$. If $f\not\in B$, then $B\subo\oplus
Ff=\pi_I(B\subo)\oplus \pi_J(B\subo)$, which is a subalgebra of
$\bigl(K_{10}\bigr)\subo$ (recall that $\bigl(K_{10}\bigr)\subo$ is
the direct sum of the ideals $I$ and $J=Ff$, here $\pi_I$ and
$\pi_J$ denote the corresponding projections), and so the subalgebra
generated by $B$ and $f$ equals $B\oplus Ff$. Since $B$ is maximal
$B\oplus Ff=K_{10}$ so, in particular,
$B\subuno=\bigl(K_{10}\bigr)\subuno$, but
$\bigl(K_{10}\bigr)\subuno$ generates $K_{10}$, so $B=K_{10}$, a
contradiction. Therefore $f\in B$, and hence $e=1-f\in B$ too.

We use now the description of $K_{10}$ given at the beginning of
Section 2 and due to D. King. Then $B_{\bar{0}} = Fe+ Ff+ V_0,$ with
$V_0$ a vector subspace of $V$ and so, in order to study the maximal
subalgebras of $K_{10}$, we analyze the possible dimensions of
$V_0.$

\bigskip

(1) If $V_0 =V$ then $B=(K_{10})_{\bar{0}}$ is a maximal
subalgebra of $K_{10}.$

\bigskip

(2) If $\dim V_0 = 1,$ we have two possible situations: either  the
rank of $Q$ is 1 or is 0.

\smallskip

If $\rank (Q \vert_{V_0}) =1$ using  Witt's Theorem  and Theorem 3.3
we can suppose that $V_0=Fa.$ Now $B_{\bar{1}}= P_B+Q_B$ with
$P_B=\{ p\in B : p a=p\}, \ Q_B=\{q\in B : q a = -q\}$. If $\dim
P_B=2$ we can check that $B= \langle e,f,a,p_1,p_2 \rangle$ is a
maximal subalgebra of $K_{10}.$ Likewise if $\dim Q_B=2,$ $B=
\langle e,f,a,q_1,q_2 \rangle$ is a maximal subalgebra of $K_{10}$
but there is an automorphism of $K_{10}$ applying $ \langle
e,f,a,p_1,p_2 \rangle$ into $\langle e,f,a,q_1,q_2 \rangle$ (the
automorphism given by the isometry which applies $a$ to $-a$ and fix
$b,$ $c_1,$ $c_2$). If $\dim P_B =\dim Q_B=1$, then we can check
that $\dim V_0> 1,$ a contradiction. If $P_B=0$ then $B\subseteq
\langle e,f,a,q_1, q_2 \rangle ,$ which is a maximal subalgebra, and
we are in the case above.

\smallskip

 If $\rank (Q \vert_{V_0})=0,$ using Witt's Theorem and Theorem 3.3
 we can suppose that $V_0=Fc_1.$ But  $c_1B_{\bar{1}} \subseteq Fp_1+Fq_1$ and $B_{\bar{1}}B_{\bar{1}} \subseteq
 Fc_1+Fe+Ff,$ so $B_{\bar{1}} \subseteq \langle p_1,q_1 \rangle.$ It follows that  $B= \langle e, f, c_1, p_1, q_1\rangle$ is not maximal since
  $B\subsetneqq \langle e,f,a,b,c_1,p_1,q_1 \rangle.$

\bigskip

(3) If $\dim V_0=2$, then  $\rank (Q \vert_{V_0})=0, 1$ or $2.$

\smallskip

 If $\rank (Q \vert_{V_0})=2$ we can suppose that $V_0= Fa+Fb,$ again by Witt's Theorem and Theorem
3.3. As in (2), we have that $B_{\bar{1}} = P_B+Q_B.$ If $\dim
P_B=2,$ then since $bP_B=Q_B,$ it follows that $\dim Q_B=2,$ and
then $B=K_{10},$ a contradiction. If $\dim P_B=1$ then $\dim Q_B=1,$
but since $B_{\bar{1}}B_{\bar{1}}\subseteq B_{\bar{0}},$ it holds
that either $B_{\bar{1}} = \langle p_1,q_2 \rangle$ or $B_{\bar{1}}
= \langle p_2,q_1 \rangle.$ This is a contradiction with
$bB_{\bar{1}}\subseteq B_{\bar{1}}.$ If $\dim P_B=0$ then $\dim
Q_B=0$ and $B$ is not a maximal subalgebra because $B\subseteq
(K_{10})_{\bar{0}}.$

\smallskip

 If $\rank (Q\vert_{V_0}) =1$ then $V_0=Fa+Fc_1$ and again $B_{\bar{1}} = P_B+Q_B.$ If $\dim
P_B=2$ then $\dim Q_B\geq 1$ (because  $c_1p_1=0$ and $c_1p_2=q_1$)
and so $B=K_{10},$ a contradiction. If $\dim P_B=1$ then $\dim
Q_B=1$ or $0.$ If $\dim Q_B=1,$ from
$B_{\bar{1}}B_{\bar{1}}\subseteq B_{\bar{0}}$ we have that
$Q_B=Fq_1$ and for $0\neq p\in P_B$ it follows that $pq_1 \in V_0=
Fa+Fc_1,$ and therefore $p\in Fp_1.$ But in this case $B$ is not
maximal since $B\subsetneqq \langle e,f,a,b,c_1,p_1,q_1 \rangle.$
And if $\dim Q_B=0,$ then $P_B=Fp_1$ and again $B$ is not maximal.
Now if $\dim P_B=0$, then $\dim Q_B=1,$ with $Q_B= \langle q_1
\rangle.$ But again $B$ is not maximal because $B\subsetneqq \langle
e,f,a,b,c_1,p_1,q_1 \rangle.$

\smallskip

 Finally, if $\rank (Q \vert_{V_0})=0,$ we can suppose that
$V_0 = F \cdot (a+b)+Fc_1,$ by Witt's Theorem and Theorem 3.3. If
$x=\alpha p_1+\beta p_2 + \sigma q_1 +\gamma q_2\in B_{\bar{1}},$
$c_1x=\beta q_1+\gamma p_1\in B_{\bar{1}}$ and $(a+b)(c_1x) =
(\gamma - \beta) (p_1- q_1) \subseteq B_{\bar{1}},$ therefore if
$\gamma \neq \beta$ then $\{p_1, q_1 \}\subseteq B_{\bar{1}}.$ So
either $\{p_1, q_1\} \subseteq B_{\bar{1}}$ or $B_{\bar{1}}
\subseteq Fp_1 + Fq_1 + F(p_2+q_2).$ Moreover with $\gamma \neq
\beta,$ $q_1x\notin B_{\bar{0}}.$ Therefore $\gamma = \beta $ and
$B= \langle e,f,a+b,c_1, p_1, q_1, p_2+q_2 \rangle,$ which is a
maximal subalgebra.

\bigskip

(4) If $\dim V_0=3$, then $\rank (Q \vert_{V_0})= 2$  or $3.$

\smallskip

 If $\rank (Q \vert_{V_0}) =3,$ then we can suppose by  Witt's
Theorem and Theorem 3.3 that $V_0 = Fa+Fb+F \cdot (c_1+c_2).$ Again
$B_{\bar{1}}= P_B+Q_B.$ If $\dim P_B=2$ then $\dim Q_B=2$ (because
$bP_B\subseteq Q_B$) and $B_{\bar{1}}=(K_{10})_{\bar{1}},$ a
contradiction. So $\dim P_B=1 = \dim Q_B,$ but then, since
$(c_1+c_2)B_{\bar{1}} \subseteq B_{\bar{1}}$ and
$B_{\bar{1}}B_{\bar{1}} \subseteq B_{\bar{0}}$, it follows that
either $B_{\bar{1}} = \langle p_2,q_1 \rangle$ or $B_{\bar{1}} =
\langle p_1,q_2 \rangle,$ and this contradicts that $bB_{\bar{1}}
\subseteq B_{\bar{1}}.$

\smallskip

 If $\rank(Q \vert_{V_0})=2$, then $V_0=Fa+Fb+Fc_1$ and $B_{\bar{1}}=P_B+Q_B,$ as always.
If $\dim P_B=2$ then $\dim Q_B=2,$ a contradiction. So $\dim P_B=1$
and $\dim Q_B=1.$ Since $P_BQ_B \subseteq V_0$ we have that
$B_{\bar{1}}= \langle p_1,q_1 \rangle$ and $B= \langle
e,f,a,b,c_1,p_1,q_1 \rangle,$ which is a maximal subalgebra of
$K_{10}.$
\end{proof}

\bigskip

\vspace{0.5cm}

Next we describe the structure of each one of the types of maximal
subalgebras that we have obtained for $K_{10}$:

\begin{enumerate}

\item[(i)] $(K_{10})_{\bar{0}}$ is semisimple, and we have seen that  $K_{10} \cong J(V, Q) \oplus F.$

\bigskip

\item[(ii)] $<e, f, a, p_1, p_2>$ is semisimple and isomorphic to $F \oplus
D_{-6}$ where \smallskip

$F \cong \langle 2e-a\rangle \thickspace \thickspace \hbox{ and }
\thickspace \thickspace \langle f, 2e+a\rangle+ \langle p_1,
p_2\rangle \cong D_{-6}.$

\bigskip

\item[(iii)] $B=<e, f, a+b, c_1, p_1, q_1, p_2+q_2>$ has the following radical
$$R = \langle c_1, a+b, p_1+q_1\rangle,$$
\noindent and the quotient of the subalgebra $B$ by its radical is
isomorphic to $D_{-3}$:
$$
B/R = \langle \bar{e}, \bar{f},
\overline{p_1}, \overline{\frac{p_2+ q_2}{2}}\rangle \cong D_{-3}.
$$

\bigskip

\item[(iv)] $B= <e, f, a, b, c_1, p_1, q_1>$ has as radical $R=
\langle c_1, p_1, q_1\rangle$ with $B/R \cong F \oplus T$ with
$T$ the Jordan algebra of a bilineal form.
$$B/R \cong \langle \bar{e}, \bar{f},
\bar{a}, \bar{b}\rangle.$$

\end{enumerate}
\bigskip

And finally we will give the description of the maximal subalgebras
of $K_{10}$ in terms of the realization shown in Section 2, given by
Benkart and Elduque. First consider the map
\begin{eqnarray*}
\varphi \colon \quad K_{10} \qquad &\longrightarrow & \qquad K_{10} \\
    \lambda 1 + a\otimes b & \longmapsto & \lambda 1 +
(-1)^{\bar a \bar b} b\otimes a,
\end{eqnarray*}
with $a, b $ homogeneous elements of $K_3.$ We remark that $\varphi$
 is an automorphism because
\begin{alignat*}{2}
\varphi (a \otimes b) \varphi (c \otimes d)  =& \thinspace ((-1)
^{\bar a \bar b}
 b\otimes a) \cdot ((-1)^{\bar c \bar d} \thinspace  d \otimes c)  \\
= & \thinspace  (-1)^{\bar a \bar b + \bar c \bar d + \bar a \bar
d} \thinspace  (bd \otimes ac -\frac{3}{4}
(b \thinspace | \thinspace d)(a \thinspace | \thinspace c) 1) \\
=& \thinspace (-1)^{\bar b \bar c}((-1)^{(\bar a + \bar c)(\bar b
+ \bar
d)}(bd \otimes ac - \frac{3}{4} (a \thinspace | \thinspace c)(b \thinspace | \thinspace d)1)) \\
= & \thinspace \varphi ((a\otimes b)(c\otimes d)).
\end{alignat*}

\smallskip

\begin{coro}
The maximal subalgebras of $K_{10}=F\cdot 1\oplus (K_3\otimes K_3)$
are, up to conjugation by automorphisms, the following:
\begin{enumerate}

\item[{\rm (i)}] $(K_{10})_{\bar{0}}.$

\item[{\rm (ii)}] The fixed subalgebra by the automorphism
$\varphi$.

\item[{\rm (iii)}] The subalgebra $F\cdot 1\oplus (M\otimes K_3)$,
where $M=F1+Fx$ is a maximal subalgebra of $K_3$.

\item[{\rm (iv)}] The subalgebra $F \cdot 1+ F(e \otimes e)+ (x \otimes K_3)+ (K_3 \otimes
y)$.
\end{enumerate}
\end{coro}
\begin{proof}

 It is easy to compute the fixed elements by $\varphi$ in terms
of the original basis of $K_{10}$ given in the introduction. This
subalgebra is  $\langle e,f,c_1+c_2, p_1-q_2,
 p_2-q_1 \rangle$ and now using Proposition 3.1 we know that there exists an automorphism
of  $K_{10}$ that applies $\langle e,f,a,p_1,p_2 \rangle$ into
$\langle e,f,c_1+c_2, p_1-q_2, p_2-q_1 \rangle$ (notice that
$F(c_1+c_2)$ and $Fa$ are isometric relative to the bilinear form
$Q$ of $V,$ and moreover $\langle p_1-q_2,p_2-q_1 \rangle = \{ x\in
(K_{10})_{\bar{0}} \ \vert \  \ (c_1+c_2)x= -x\}$). So the maximal
subalgebras of type (ii) in Theorem 4.1 are conjugated to the
subalgebra of fixed elements by $\varphi$.

\smallskip

Also the maximal subalgebras of type (iii) correspond, through the
isomorphism in Theorem \ref{th:BeEl}, to $\langle 1, e\otimes e,
x\otimes x, x\otimes y, e\otimes x, x\otimes e, e\otimes y \rangle,$
that is, to subalgebras $F \cdot 1 \oplus (M \otimes K_3)$ with $M$
maximal subalgebra of $K_3$; while the maximal subalgebras of type
(iv) correspond, through the isomorphism in Theorem \ref{th:BeEl},
to
 $F \cdot 1+ F(e \otimes e)+ (x \otimes K_3)+ (K_3 \otimes y).$
\end{proof}

\end{document}